\documentclass[11pt]{amsart}
\usepackage{amssymb,color}
\usepackage{graphicx}
\usepackage{amssymb}
\usepackage{epstopdf}
 \usepackage{latexsym}
\usepackage{amsmath}
\usepackage{amsfonts} 
\usepackage{amsthm}

\newtheorem{Proposition}{Proposition}
\newtheorem{Theorem}[Proposition]{Theorem}

\def\XXint#1#2#3{{\setbox0=\hbox{$#1{#2#3}{\int}$}
\vcenter{\hbox{$#2#3$}}\kern-.5\wd0}}

\def\a{\alpha}

 \def\Box{{\hfill\hbox{\enspace${\sqre}$}} \smallskip}
    \def\sqr#1#2{{\vcenter{\vbox{\hrule height .#2pt
                             \hbox{\vrule width .#2pt height#1pt \kern#1pt
                                   \vrule width .#2pt}
                             \hrule height .#2pt}}}}
 \def\sqre{\mathchoice\sqr54\sqr54\sqr{4.1}3\sqr{3.5}3}

     \def\CC{\mathbb{C}}
    
    \def\NN{\mathbb{N}}

    \def\ZZ{\mathbb{Z}}

\def\be{\begin{equation}}
\def\ee{\end{equation}}

\def\al{\alpha}
\def\nint{\not\negmedspace{\negthinspace\int}}

\def\b{\beta} 
\def\a{\alpha}
\def\be{\begin{equation}}
\def\bel{\begin{equation}\label}
\def\ee{\end{equation}}

 \def\CC{\mathbb{C}}

\def\ga{c}

\def\XXint#1#2#3{{\setbox0=\hbox{$#1{#2#3}{\int}$ }
\vcenter{\hbox{$#2#3$ }}\kern-.6\wd0}}

\begin{document} 

%\today
%\date{ }
\title[Jacobi series for general parameters ]{Jacobi series for general parameters and applications}
\author{Rodica D. Costin, Marina David}
\address{Department of Mathematics\\The Ohio State University\\Columbus, OH 43210}

\maketitle

\begin{abstract}

The representation of analytic functions as convergent series in Jacobi polynomials $P_n^{(\a,\b)}$ is reformulated using a unified approach for all $\a,\b\in\CC\setminus\{0,-1,-2,\ldots\},\ \a+\b\ne -2,-3,\ldots$ The coefficients of the series are given as usual integrals in the classical case (when $\Re\al,\Re\beta>-1$), or by the Hadamard principal part of these integrals when they diverge.

As an application it is shown that nonhomogeneous differential equations of hypergeometric type do generically have a unique solution which is analytic at both singular points in $\CC$.

\end{abstract}

%Subject classification: Primary 33C45, 42C10 Secondary: 35A25

\section{Introduction}

The purpose of this note is to provide a unified and easy-to-use formalism of expansion of analytic functions in Jacobi series for general parameters. 
Such expansions are useful in approximation theory, as well as in a variety of other areas, such as the study of differential equations (see e.g. \cite{RDCMDK}, \cite{RDC0}, \cite{CorLin}) and non-selfadjoint spectral problems which appear in the study of stability of nonlinear partial differential equations.  An application answering another question on differential equations is given in \S\ref{hyper}.

\subsection{An overview of Jacobi polynomials for classical parameters} We summarize below a few well known facts about Jacobi polynomials $P_n^{(\a,\b)}(z),\ n=0,1,2,\ldots$ for $\a,\b>-1$.  (See \cite{Ismail} or \cite{nist} for an overview). 

The Jacobi polynomials are determined by the condition of mutual orthogonality on the interval $[-1,1]$ with respect to the weight 
$$W(z)=(1-z)^\a(1+z)^\b$$
namely
\bel{defJacobi}
\int_{-1}^1P_n^{(\a,\b)}(z)\, P_k^{(\a,\b)}(z)\, W(z)\,dz=\mathcal{A}_n\,\delta_{n,k}
\ee
where
\bel{valAn}
\mathcal{A}_n=2^{\a+\b+1}\frac{\Gamma(n+\a+1)\, \Gamma(n+\b+1)}{ (2n+\a+\b+1)\,\Gamma(n+\a+\b+1)\, n!}\ \ \text{for }n=0,1,2,\ldots
\ee
Relation \eqref{defJacobi} and the condition that the coefficient of the leading term, $z^n$ in $P_n^{(\a,\b)}(z)$ be positive uniquely determine the Jacobi polynomials.

In the particular case when $\a=\b$ these polynomials are called Gegenbauer (or {\em ultraspherical}), and, when $\a,\b=\pm\frac12$ they are called Chebyshev polynomials.

As any orthogonal polynomials, Jacobi polynomials satisfy a three term recurrence relation
\bel{ttr}
P_{n+1}^{(\a,\b)}(z)=(A_nz+B_n) P_n^{(\a,\b)}(z)-C_nP_{n-1}^{(\a,\b)}(z)
\ee
where $P_0^{(\a,\b)}(z)=1,\, P_1^{(\a,\b)}(z)=A_0z+B_0$ and
$$A_n=\frac{(2n+\a+\b+1)(2n+\a+\b+2)}{2(n+1)(n+\a+\b+1)},$$

$$B_n=\frac{(\a^2-\b^2)(2n+\a+\b+1)}{2(n+1)(n+\a+\b+1)(2n+\a+\b)},$$

$$C_n=\frac{(n+\a)(n+\b)(2n+\a+\b+2)}{(n+1)(n+\a+\b+1)(2n+\a+\b)}.$$

Also, the Jacobi polynomials satisfy the Rodrigues formula:
\bel{Rodr}
 P_n^{(\alpha,\beta)}(z)= D_n\, \frac{1}{W(z)}\, \frac{d^n}{dz^n}\left[ Q(z)^nW(z)\right]
 \ee
where 
\bel{WQC}
W(z)=(1-z)^\a(1+z)^\beta,\ Q(z)=(1-z)(1+z),\ \ D_n=\frac{(-1)^n}{2^nn!}.
\ee

Furthermore, $P^{(\a,\b)}_{n}(z)$ satisfy the differential equation:
\begin{equation}\label{lamJacobi}
(1-z^2)y''+\left[ \b-\a-(\a+\b+2)z\right])y'+n(n+\a+\b+1)\,y=0
\end{equation}
i.e. $P^{(\a,\b)}_{n}(z)$ are eigenfunctions of a self-adjoint operator.

It is a classic result that analytic functions can be expanded in Jacobi series (see \cite{Szego} p.245):

\begin{Theorem}\label{SZG}
 Let $\a>-1,\ \beta>-1$. Let $f$ be analytic in the interior of an ellipse with foci at $\pm1$. Denote by $\mathcal{E}$ the greatest such ellipse.
Then $f$ has an expansion 
\begin{equation}\label{serf1}
f(z)=\sum_{n=0}^\infty f_n P_n^{(\alpha,\beta)}(z)
\end{equation}
which is convergent in the interior of $\mathcal{E}$ and divergent outside.
\end{Theorem}

Of course, due to the orthogonality relation \eqref{defJacobi} the coefficients $f_n$ in \eqref{serf1} can be calculated as

\bel{ak}
f_n=\mathcal{A}_n^{-1}\, \int_{-1}^1 \, f(z)\,  P_n^{(\alpha,\beta)}(z)\,W(z)\, dz
\ee

\subsection{Jacobi polynomials for general parameters} 
For $\a,\b\in\CC$ with $\Re\a,\Re\b>-1$ all the results presented above still hold - except that now $P^{(\a,\b)}_{n}(z)$ are eigenfunctions of a non-selfadjoint operator.

For other complex values of $\a,\b$ however, the integrals in \eqref{defJacobi} and in \eqref{ak} diverge and the Jacobi polynomials can no longer be thought of being an orthogonal system. 
 Nevertheless, Jacobi polynomials can be defined through the Rodrigues formula, or, most importantly, through the three-terms recurrence. The latter then implies, by Favard's Theorem, that there is a bilinear form with respect to which the polynomials are orthogonal \cite{Favard}, \cite{Sho}.

Many decades after Theorem\,\ref{SZG} has been established, in his 1974 paper \cite{Carlson}, Carlson showed that expansions \eqref{serf1} hold for general complex parameters $\a,\b$ provided that $\a+\b\ne -2,-3,-4,\ldots$. The formulas in \cite{Carlson} use an integral kernel and are not easy to use in practice (see Theorem\,\ref{ThCarl} in \S\ref{ThC}).
More recently Kuijlaars, Mart\'inez-Finkelshtein and Orive find by analytic continuation that
orthogonality of Jacobi polynomials can be established by integration on special paths in the complex plane; they also derive an associated Riemann-Hilbert problem; in some cases incomplete or quasi-orthogonality are found, or even multiple orthogonality conditions \cite{9}, \cite{10}.

Later, in \cite{RDC}, it was shown that for $\a,\b\in\CC\setminus\ZZ_-$ with $ \a+\b\ne -2,-3,-4\ldots$ the Jacobi polynomials do satisfy \eqref{defJacobi} only now the integral should be understood as its Hadamard finite part: \eqref{defJacobi} taken in the Hadamard finite part is the bilinear form with respect to which the Jacobi polynomials are orthogonal.  It was also shown in \cite{RDC} that the Hadamard finite part of the integral is the analytic continuation in $\a,\b$ of the bilinear form \eqref{defJacobi} for the cases when the integral no longer converges.

\subsection{Jacobi series for general parameters through the Hadamard finite part}
We reformulate both the classical Theorem\,\ref{SZG} and Carlson's result from a unified point of view, easy to use in applications: Theorem\,\ref{propAp} states that even in the case of general complex parameters \eqref{defJacobi} and \eqref{ak} are valid provided that possible divergent integrals are replaced by their Hadamard finite part. A brief introduction to the Hadamard finite part  is found in the Appendix.

The result is formulated for Jacobi polynomials on the interval $[0,1]$:

\begin{Theorem}\label{propAp}

Let $\a,\beta\in \CC\setminus\{ 0,-1,-2,\ldots\},\ \a+\b\ne -2,-3,\ldots$

Consider the (non-normalized) Jacobi polynomials on $[0,1]$ given by the Rodrigues formula:
\begin{equation}\label{Rodrpn}
p_n=\frac{1}w \, \frac{d^n}{dx^n}\left(q^nw\right)\ \ \text{where }w(x)=x^\a(1-x)^\beta,\ q(x)=x(1-x).
\end{equation}

(i) {\em Orthogonality:}
The polynomials $p_n$are orthogonal with respect to a bilinear form:
\bel{ortogonnn}
\mathcal{H}\!\!\!\int_0^1 p_n(x)\,p_k(x)\, w(x)\, dx=\delta_{nk}\, {a}_n
\ee
where $\mathcal{H}\!\! \int$ denotes the Hadamard finite part of the integral and
\bel{consta}
a_n=\frac{(n!)^2}{2^{\a+\b+1}}\,\mathcal{A}_n=\frac{ n!\, \Gamma(n+\a+1)\, \Gamma(n+\b+1)}{ (2n+\a+\b+1)\,\Gamma(n+\a+\b+1)}.
\ee

(ii) {\em Completeness:}
Let $\Omega\in\CC$ be the interior of an ellipse with foci $1$ and $0$, and let $f$ be analytic in $\Omega$.  Then $f(x)$ has the expansion 
\begin{equation}\label{serfgen}
f(x)=\sum_{n=0}^\infty f_n\, p_n(x),\ \ \ x\in\Omega
\end{equation}
where
\bel{formfn}
 f_n=\ a_n^{-1}\ \mathcal{H}\!\!{\int_0^1fp_nw}.
 \ee

The expansion \eqref{serfgen} converges absolutely and uniformly on every compact set in $\Omega$.
\end{Theorem}

The proof of Theorem\,\ref{propAp} is contained in \S\ref{Proof}. This result is used in \S\ref{hyper} to investigate existence of eigenfunctions for some non-self-adjoint operators; more precisely, we establish existence of unique analytic solutions of inhomogeneous differential equations of hypergeometric type under precise conditions on the parameters.

%%%%%%%%%%%%%%%%%%%%%%%%%%%%
\section{Proof of Theorem\,\ref{propAp}}\label{Proof}
\subsection{Jacobi polynomials on $[0,1]$}
Substituting 
\bel{pnPn}
p_n(x)=D_n^{-1}\, (-2)^{-n}\,P_n^{(\alpha,\beta)}(1-2x)
\ee
 in the Rodrigues formula \eqref{Rodr} we obtain \eqref{Rodrpn}
which shows that the polynomials defined by \eqref{Rodrpn} are indeed Jacobi polynomials (up to multiplicative constants).

It is easy to see that the leading coefficient of $p_n$ is the coefficient of $x^n$ which equals
$$k_n=(-1)^n(n+\a+\b+1)_n$$
(Note that if $\a+\b$ is a negative integer these polynomials are no longer a complete system).

It is also easy to see that $p_n$ depend polynomially on $\a$ and $\b$.

\subsection{Orthogonality for general parameters} Formula \eqref{ortogonnn} was proved in \cite{RDC}. 
 But one can also reason that it follows by analytic continuation, in the following way. For $\Re\a,\Re\b>-1$ substituting \eqref{pnPn} in \eqref{defJacobi} we obtain
\begin{equation}\label{ortop}
\int_{0}^1\, p_n(x)\, p_k(x)\, w(x)\, dx=0\ \ \text{if }n\ne k,
\end{equation}
a relation which, by analytic continuation, then holds for all other complex values of $\a,\b$, implying \eqref{ortogonnn}.

Formula \eqref{consta} can be obtained by a straightforward calculation based on \eqref{pnPn} and \eqref{defJacobi}, \eqref{valAn}.

\subsection{Completeness}\label{ThC} Since it was established in \cite{RDC} that analytic continuation in $\a,\b$ of integrals of the type \eqref{defJacobi}, \eqref{ak}  is given by their Hadamard finite part, then Theorem\,\ref{propAp} (ii) follows by analytic continuation of Theorem\,\ref{SZG}.

We will however give a careful proof based on the result in \cite{Carlson}, which is the following:   
\begin{Theorem}\label{ThCarl} (Theorem 1.1 in \cite{Carlson} for  $r=1,\, s=0$)
Let $\Omega$ be an open elliptic disk with foci $1$ and $0$, and let $f$ be analytic on $\Omega$. Let $(\a,\b)\in\CC^2$ and assume $\a+\beta\ne -2,-3,-4,\ldots$. Then, for every $z\in\Omega$,
\begin{equation}\label{serf}
f(z)=\sum_{n=0}^\infty\frac{1}{n!}\, c_n\, R_n(z)
\end{equation}
where  $R_n(z)$ are scalar multiples of $P_n^{(\a,\beta)}(\frac{1-z}{2})$, given by the formula $R_n(z)=R_n(-\a-n,-\beta -n,z-1,z)$ where
\begin{equation}\label{defRn}
R_n(b,b',x,y)=\frac 1{B(b,b')}\,\int_0^1\left[ux+(1-u)y\right]^n\, u^{b-1}(1-u)^{b'-1}\, du
\end{equation}
for $\Re b,\Re b'>0$ and its analytic continuation for other complex values of $b,b'$, and 
$c_n=F^{(n)}(1+\a+n,1+\beta+n;1,0)$ where 
\be
F^{(n)}(b,b';1,0)=n!(2\pi i)^{-1}\int_\gamma f(w)R_{-n-1}(b,b';w-1,w)\, dw
\ee
where $\gamma$ is a rectifiable Jordan curve in $\Omega$ which encircles the segment $[0,1]$ in the positive direction.

The series \eqref{serf} converges absolutely on $\Omega$, uniformly on every compact set in $\Omega$.

\end{Theorem}

%As always, $B(b,b')$ denotes the beta function: $$B(b,b')=\int_0^1t^{b-1}(1-t)^{b'-1}\, dt=\frac{\Gamma(b)\,\Gamma(b')}{\Gamma(b+b')}$$

A straightforward calculation can be also made to link \eqref{serfgen} and \eqref{serf}, yielding
$$p_n(x)= (-1)^n\, (\a+\beta+2n)_n\, R_n(x)$$
and
$$ \frac 1{n!}\, c_n\, (-1)^n\, (\a+\beta+2n)_n^{-1}=\left({ \mathcal{H}\!\!\!\int_0^1p_n^2w}\right)^{-1}\,{\mathcal{H}\!\!\!\int_0^1fp_nw}\ :=f_n$$
where $R_n(x)$ are given by \eqref{defRn}, $(a)_n$ denotes the rising factorial: $(a)_n=a(a+1)...(a+n-1)=\Gamma(a+n)/\Gamma(a)$. We illustrate below the main steps of the calculation.

In \cite{Carlson} $R_n$ is shown to be a multiple of a Jacobi polynomial by using the binomial formula in \eqref{defRn}, then integrating term by term (this is also correct if one uses the Hadamard principal part in case the integrals diverge). Then one compares with the expression of $p_n$ obtained by using Leibniz's rule to expand  its Rodrigues formula.

Integration by parts also holds for the Hadamard finite part, and $p_n$ were proved to be an orthogonal set in \cite{RDC}. It only remains to calculate the "norm" of $p_n$ in this sense.

We first calculate the dominant term in $p_n$: it is found by retaining only the highest powers of $x$ in its Rodrigues formula, yielding:
$p_n=D_nx^n+...$ where $D_n=(-1)^n(\a+\beta+2n)_n$.

Then, using the fact that $p_n$ is orthogonal to all the polynomials of degree less than $n$, and then using the Rodrigues formula and integrating by parts $n$ times we obtain
\begin{multline}
\mathcal{H}\!\!\!\int_0^1p_n^2w=D_n\, \mathcal{H}\!\!\!\int_0^1p_n\, x^n\, w=D_n\, \mathcal{H}\!\!\!\int_0^1\left[q^nw\right]^{(n)}\, x^n\\
=D_n\, (-1)^n\, n! \, \mathcal{H}\!\!\!\int_0^1x^{n+\a}(1-x)^{n+\beta}\, dx=D_n\, (-1)^n\, n!\, B(n+\a+1,n+\beta+1)\hskip 2cm
\end{multline}
$\Box$\\

We will need the differential equation \eqref{lamJacobi} transcribed for the Jacobi polynomials on $[0,1]$:

\begin{Proposition}\label{geneif}
The Jacobi polynomials $p_n$ defined by \eqref{Rodrpn} for general $\a,\b$ satisfy the differential equation
\begin{equation}\label{deqpn}
y''+\bigg(\frac{\a+1}{x}+\frac{\b+1}{x-1}\bigg)\,y'+\frac{n(n+\a+\b+1)}{x(1-x)}\,y=0
\end{equation}

\end{Proposition}

The proof is a straightforward calculation using \eqref{pnPn} in  \eqref{lamJacobi}. $\Box$

%%%%%%%%%%%%%%%%%%%%%%%%%%%%
\section{Analytic solutions of inhomogeneous hypergeometric equations}\label{hyper}

The differential equation of hypergeometric type
\begin{equation}\label{eqhyp} 
w''(x)+\Bigg(\frac{a}{x}+\frac{b}{x-1}\Bigg) w'(x)+\frac{c}{x(1-x)}w(x)=0
\end{equation}
has solutions which are analytic at the singular point $x=0$ and solutions analytic at $x=1$ (as it is well known, and easy to establish using Frobenius theory \cite{Codd}). But, for generic parameters $a,b,c$ no non-zero solution is analytic at both singular points $x=0$ and $x=1$.

We investigate if an inhomogeneous version, eq. \eqref{eeqq} below, has solutions that are analytic at both $x=0$ and $x=1$. This type of questions appear as generalized eigenfunction problems in many types of applications, for instance in the study of blow-up of wave maps \cite{Bison}, \cite{OC}.
 
We find the answer to be positive:
 \begin{Theorem}\label{T1}

Consider the following differential equation of hypergeometric type with an nonhomogeneous term:
\begin{equation}\label{eeqq}
u''(x)+\left(\frac{a}{x}+\frac{b}{x-1}\right) u'(x)+\frac{c}{x(1-x)}u(x)=\frac{g(x)}{x(1-x)}
\end{equation}
where $g$ is a function that is analytic in the domain $\mathcal{E}$ which is the interior of an ellipse with foci at $0,1$.

Assume that: $a,b\in\CC\setminus\{1,0,-1,-2,-3,\ldots\},\ a+b\ne 0,-1,-2,\ldots$ and
\begin{equation}\label{solcon}
c \neq n(n+a+b-1)\ \ \ \text{for all }n=0,1,2,3,\ldots
\end{equation}

Then equation \eqref{eeqq} has a unique solution which is analytic at both $x=0$ and $x=1$, and therefore it is analytic on the whole domain where $g$ is analytic.

\end{Theorem}

{\em Proof.}

Let $p_n$ be the Jacobi polynomials on $[0,1]$ defined by \eqref{Rodrpn} with $\a=a-1,\,\b=b-1$.

Denote
\begin{equation}\label{notlamn}
 \lambda_n=n(n+a+b-1).
\end{equation}

Then by \eqref{deqpn} each polynomial $p_n$ satisfies
\begin{equation}\label{pneigenfun}
p_n''(x)+\bigg(\frac{a}{x}+\frac{b}{x-1}\bigg)\,p_n'(x)=-\frac{\lambda_n}{x(1-x)}\,p_n(x)
\end{equation}
which implies that
\begin{equation}\label{pnide}
p_n''(x)+\bigg(\frac{a}{x}+\frac{b}{x-1}\bigg)\,p_n'(x)+\frac{\ga}{x(1-x)}\, p_n(x)=\frac{\ga-\lambda_n}{x(1-x)}\,p_n(x)
\end{equation}
so that $p_n(x)$ is a solution of \eqref{eeqq} for the case when $g(x)=(\ga-\lambda_n)p_n(x)$.

We are thus led to the idea of expanding the inhomogeneous term $g$ in a series in Jacobi polynomials.

By Theorem\,\ref{propAp} there exist constants $g_0,g_1,g_2,\ldots$ so that $
g(x)= \sum_{n=0}^\infty g_n\, p_n(x)$ with the series converging uniformly on every compact contained in $\mathcal{E}$ and $g_n$ given by the formula
\bel{formgn}
 g_n=a_n^{-1}\ \mathcal{H}\!\!{\int_0^1gp_nw}
 \ee

We look for a solution $u(x)$ of \eqref{eeqq} which is analytic at both $0$ and $1$, and therefore it is analytic in $\mathcal{E}$; also by Theorem\,\ref{propAp}, such a solution has an expansion
$u(x)= \sum_{n=0}^{\infty} u_n\, p_n(x)$ which plugged in the equation \eqref{eeqq} yields
\begin{equation}
\sum_{n=0}^\infty u_n \left[\, p_n''+\left(\frac{a}{x}+\frac{b}{x-1}\right)\, p_n'+\frac{\ga}{x(1-x)}\,p_n\right]=\sum_{n=0}^{\infty} \frac{g_n}{x(1-x)} \, p_n(x)
\end{equation}
which by \eqref{pnide} implies 
$$\sum_{n=0}^\infty\, \frac{\ga-\lambda_n}{x(1-x)}\,u_n\, p_n(x)=\sum_{n=0}^{\infty} \frac{g_n}{x(1-x)} \, p_n(x)
$$
and using \eqref{ortogonnn} it follows that $u_n\, (c-\lambda_n)=g_n$ for all $n$. In view of \eqref{solcon} and \eqref{notlamn} it follows that 
\begin{equation}\label{uk}
u_n = \frac{g_n}{\ga-n(n+a+b-1)},\ \ \ n=0,1,2,\ldots
\end{equation}

We thus found a unique solution analytic in $\mathcal{E}$ of \eqref{eeqq} as the sum of the Jacobi series

\begin{equation}
u(x)=\sum_{n=0}^\infty\frac{g_n}{\ga-n(n+a+b-1)}\, p_n(x)
\end{equation}
with $g_n$ given by \eqref{formgn}. The series is obviously convergent everywhere the series of $g$ is. $\Box$

%%%%%%%%%%%%%%%%%%%%%%%%%%%%%%%

\section{Appendix: The Hadamard finite part}\label{Append}

The concept of the {\em{finite part}} of a (possibly divergent) integral was introduced by Hadamard \cite{Hadamard} as a convenient way to express solutions of differential equations. He showed that this finite part of an integral (which coincides with the usual value if the integral is convergent) can be combined and manipulated in much the same way as usual integrals: they are additive on the interval of integration, changes of variable are allowed, etc. (They do not behave well with respect to inequalities.) The finite part can be calculated either by Taylor series, or by integration along closed paths in the complex plane.

Subsequently the Hadamard finite part has been interpreted in terms of distributions (see, e.g. \cite{Morton-Krall}) and it turned out that many problems of mathematical physics have solutions expressible as the Hadamard finite part of (divergent) integrals, and numerical methods of calculations have been subsequently developed (see for example \cite{Davis_Rabinowitz}).

The present section contains some properties of the Hadamard finite part of integrals of the type $\int_0^xt^{\alpha-1}f(t)\, dt$ with $f$ analytic at $0$; when $\Re\alpha\leq 0$, $\alpha\not\in(-\NN)$, its Hadamard finite part is denoted here by
$$\mathcal{H}\!\!\!\int_0^{x}t^{\alpha-1}f(t)\, dt$$

Consider $f(x)$ a function analytic at $x=0$: $f(x)=\sum_{n=0}^\infty c_n x^n$ a series with nonzero radius of convergence $R$. The integral
\begin{equation}\label{intdef}
x^{-\a}  \int_0^x t^{\a-1}f(t)\, dt
\end{equation}
converges for $\Re \a>0$ and for $|x|<R$ we have
\begin{multline}\label{calc}
x^{-\a}  \int_0^x t^{\a-1}f(t)\, dt=x^{-\a}  \int_0^x t^{\a-1} \, \sum_{n=0}^\infty c_n t^n\, dt=x^{-\a} \sum_{n=0}^\infty c_n \int_0^x t^{n+\a-1}\, dt \\
=x^{-\a} \sum_{n=0}^\infty \frac{c_n}{n+\a}\, x^{n+\a}= \sum_{n=0}^\infty \frac{c_n}{n+\a}\, x^{n}
\end{multline}
and we see that the final result is analytic at $x=0$. For all other complex values of $\a\ne 0,-1,-2,\ldots$, the Hadamard finite part of the integral is defined by:
$$x^{-\a}  \mathcal{H}\!\!\!\int_0^x t^{\a-1}f(t)\, dt\,:=\,\sum_{n=0}^\infty \frac{c_n}{n+\a}\, x^{n}\ \ \ \text{for all }\a\in\CC\setminus\{ 0,-1,-2,\ldots\}$$

This definition represents the analytic continuation in $\a$ of the function \eqref{intdef} to $\a$ in the complex plane ;when $\a$ is a negative integer the continuation has poles of order one.

{\em Remark.} For practical calculations one does not need to expand $f$ in its full power series, just a finite number of terms are needed. For example, to calculate
$$x^{1/2}\mathcal{H}\!\!\!\int_0^x t^{-3/2}f(t) dt$$
first write $f(t)=f(0)+t\tilde{f}(t)$ where $\tilde{f}$ is analytic at $0$. Then
$$x^{1/2}\mathcal{H}\!\!\!\int_0^x t^{-3/2}f(t) dt=x^{1/2}\mathcal{H}\!\!\!\int_0^x t^{-3/2}f(0) dt +x^{1/2}\int_0^x t^{-1/2}\tilde{f}(t) dt$$
where the first term is
$$x^{1/2}\mathcal{H}\!\!\!\int_0^x t^{-3/2}f(0) dt =x^{1/2}\,\frac{x^{-1/2}}{-1/2}f(0) =-2f(0) $$
and the last term is an usual integral.

%{\em II. The Hadamard finite part of integrals $\int_0^1t^{\alpha-1}(1-t)^{\b-1}f(t)\, dt$}
%Consider $f$ analytic on $[0,1]$. Since the Hadamard finite part is additive over the interval of integration \cite{RDC} we can break the integral 


\begin{thebibliography}{99}

%\bibitem{Anosov-Arnold} D. V. Anosov, V. I. Arnold (Eds.), {\em{Dynamical Systems I}} Springer-Verlag, 1988 





% 10-21

\bibitem{Bison} Biz\'on, Piotr. An unusual eigenvalue problem. {\em Acta Phys. Polon. B} {\bf 36} (2005), no. 1, 5--15.  MR2125332

\bibitem{Carlson} Carlson, B. C. Expansion of analytic functions in Jacobi series. {\em SIAM J. Math. Anal.} {\bf 5} (1974), 797--808.  MR0387692

\bibitem{Codd}  Coddington, Earl A.; Levinson, Norman
Theory of ordinary differential equations. {\em McGraw-Hill Book Company, Inc., New York-Toronto-London}, 1955. xii+429 pp.MR0069338

\bibitem{OC} Costin, O.; Donninger, R.; Xia, X. A proof for the mode stability of a self-similar wave map. {\em Nonlinearity} {\bf 29} (2016), no. 8, 2451--2473.  MR3538419

\bibitem{RDCMDK} Costin, Rodica D.; Kruskal, Martin D. Nonintegrability criteria for a class of differential equations with two regular singular points. {\em Nonlinearity} {\bf 16 }(2003), no. 4, 1295--1317.  MR1986296

\bibitem{RDC0}  Costin, Rodica D. Analytic linearization of nonlinear perturbations of Fuchsian systems. {\em Nonlinearity} {\bf 21} (2008), no. 9, 2083--2097. MR2430662 

\bibitem{CorLin} Costin, Rodica D.  Nonlinear perturbations of Fuchsian systems: corrections and linearization, normal forms. {\em Nonlinearity} {\bf 21}, (2008), no.9, 2073-2082. MR2430661

\bibitem{RDC}  Costin, Rodica D. Orthogonality of Jacobi and Laguerre polynomials for general parameters via the Hadamard finite part. {\em J. Approx. Theory} {\bf 162} (2010), no. 1, 141--152. MR2565830 








\bibitem{Davis_Rabinowitz} Davis, Philip J.; Rabinowitz, Philip. Methods of numerical integration. Second edition. Computer Science and Applied Mathematics. {\em Academic Press, Inc., Orlando, FL,} 1984. {\rm xiv}+612 pp. ISBN: 0-12-206360-0 MR0760629

\bibitem{Favard} Favard, J. Sur les polynomes de Tchebicheff. {\em C. R. Acad. Sci. Paris} (in French), {\bf 200}  (1935) 2052-2053, JFM 61.0288.01


\bibitem{Hadamard} Hadamard, Jaques. Lectures on Cauchy's problem in linear partial differential equations. {\em Dover Publications, New York,} 1953. {\rm iv}+316 pp. 36.0X MR0051411

\bibitem{Ismail} Ismail, Mourad E. H. Classical and quantum orthogonal polynomials in one variable. With two chapters by Walter Van Assche. With a foreword by Richard A. Askey. Encyclopedia of Mathematics and its Applications, 98. {\em Cambridge University Press, Cambridge}, 2005. {\rm xviii}+706 pp. ISBN: 978-0-521-78201-2; 0-521-78201-5. MR2191786


\bibitem{Hadamard2} Krommer, Arnold R.; Ueberhuber, Christoph W. Computational integration. {\em Society for Industrial and Applied Mathematics (SIAM), Philadelphia, PA}, 1998. xx+445 pp. ISBN: 0-89871-374-9 MR1625683 

\bibitem{9} Kuijlaars, A. B. J.; Martinez-Finkelshtein, A.; Orive, R. Orthogonality of Jacobi polynomials with general parameters. {\em Electron. Trans. Numer. Anal.} {\bf 19} (2005), 1?17. MR2149265


\bibitem{10} Mart\'inez-Finkelshtein, A.; Orive, R. Riemann-Hilbert analysis of Jacobi polynomials orthogonal on a single contour. {\em J. Approx. Theory} {\bf 134} (2005), no. 2, 137-170. MR2142296


\bibitem{Morton-Krall} Morton, Robert D.; Krall, Allan M. Distributional weight functions for orthogonal polynomials. {\em SIAM J. Math. Anal.} {\bf 9} (1978), no. 4, 604--626. MR0493141


\bibitem{Sho} Shohat, J.  Sur les polyn\^omes orthogonaux g\'en\'ralises. {\em C. R. Acad. Sci. Paris} (in French), {\bf 207} (1938) 556-558, Zbl 0019.40503


\bibitem{Szego} Szeg\"o, Gabor. Orthogonal Polynomials. American Mathematical Society Colloquium Publications, v. 23. {\em American Mathematical Society, New York,} 1939 {\rm ix}+401 pp. MR0000077


\bibitem{nist} Digital Library of Mathematical Functions, http://dlmf.nist.gov/18

\end{thebibliography}
\end{document}